\documentclass{article}
\usepackage{graphicx} 
\graphicspath{{./figures/}}

\usepackage{amsmath}
\usepackage{amssymb}
\usepackage{amsthm}
\usepackage{algorithm}
\usepackage{algpseudocode}

\newcommand{\gnorm}[1]{\left|\left|#1\right|\right|}
\newcommand{\cf}{{\it cf.\ }}

\title{Weak-form Extended Dynamic Mode Decomposition}
\author{Christopher W. Curtis \& David M. Bortz}
\date{}
\begin{document}

\maketitle

\begin{abstract}
In this work, we develop a weak-form version of Extended Dynamic Mode Decomposition that we call WEDMD.  We establish a number of analytic results about the method and show explicitly how the weak form is able to mitigate the impacts of noise in linear stochastic differential equations.  In nonlinear systems, we likewise show how the method is able to generate excellent approximations to Koopman modes without making recourse to sparsity promoting approaches.  This excellent approximation property allows for meaningful forecasting even in the presence of noisy data.    
\end{abstract}

\section{Introduction}

With the contemporary shift in focus towards data-driven methodologies throughout the sciences and engineering, for those problems in which data is generated by some underlying dynamical system, two broad classes of methods have seen extensive development and become ubiquitous in use across a wide range of applications.  The first class consists of generating time-evolving models without approximating the underlying vector field.  This is done by way of approximating the affiliated Koopman \cite{LasotaMackey1994ChaosFractalsNoise, BruntonBudisicKaiserKutz2022ModernKoopman, Mezic2021KoopmanOperatorGeometryLearning} operator and gives rise to a range of methods generally known as Dynamic-Mode Decomposition (DMD) \cite{Schmid2009DMD, williams, williams2, KutzBruntonBruntonProctor2016DMD}.  While primarily used for modal analysis \cite{KutzBruntonBruntonProctor2016DMD} akin to the use of principal orthogonal decompositions (PODs) in fluid mechanics and related fields \cite{berkooz, TowneSchmidtColonius2018SPOD}, advances in machine learning have allowed for DMD to be used in a time-stepping fashion for interpolation and prediction \cite{LuschBruntonKutz2018DeepLearningKoopman, AzencotErichsonLinMahoney2020ConsistentKoopmanAutoencoders, curtis_dldmd, curtis_hldmd}.  

The second class of data-driven methods in dynamical systems can be described as model-discovery algorithms that generate direct approximations to the underlying vector field.  While a broad field, in this work, we focus on dictionary-based methods that are generally referred to via Sparse Identification of Nonlinear Dynamics (SINDy) \cite{BruntonProctorKutz2016SINDy}.  The original SINDy work has seen a number of developments \cite{KahemanBruntonKutz2020PISINDy, KahemanKutzBrunton2021SINDyPI} that have helped improve stability, reliability, and account for underlying symmetries due to underlying physical constraints in the data.  However, as seen in \cite{MessengerBortz2021WSINDy, MessengerBortz2021_WSINDyPDE}, when trying to account for noise, it is more advantageous to move to a weak-formulation of the model discovery problem that is otherwise attacked directly in SINDy methods.  This weak form approach has given rise to WSINDy and a host of related methods \cite{BortzMessengerDukic2023_WENDy} that have proven to be exceptionally powerful and extendable.  

In this work, with an eye towards better management of noise, we ask how one can generate approximations of the Koopman operator from a weak viewpoint.  Proceeding from the assumptions used in extended DMD (EDMD), \cf \cite{williams, williams2}, and then applying the weak-form or Galerkin approximation approach found in \cite{MessengerBortz2021WSINDy}, we generate an algorithm that we call Weak-form Extended Dynamic Mode Decomposition (WEDMD).  This complements a growing body of results exploring weak formulations of Koopman frameworks of dynamical systems; see \cite{zhou2025koopman, xu2025datadriven, bennett2026weakdmd,colbrook2026prone}.

Two recent works warrant a sharper comparison, since both share
vocabulary with ours while addressing otherwise different problems.  In \cite{bennett2026weakdmd}, the authors develop what they call the {\it weak} DMD (WDMD) method by applying a temporal weak form (test and trial bases with integration-by-parts) to the linear problem $\dot{\mathbf{y}} = A\mathbf{y}$.  They then estimate the state by a trial-basis expansion thereby passing to an otherwise now standard DMD method used to recover the eigenvalues of the state operator $A$.  By focusing instead on building approximations to an underlying Koopman operator via a generic dictionary of observables, WEDMD allows for far greater flexibility and accuracy.  Moreover, our approach also makes the error analysis more straightforward, and this provides meaningful insights as we show.  


Likewise, the recent PRONE framework~\cite{colbrook2026prone} unifies DMD, EDMD, SINDy, and Koopman regression as a single Petrov--Galerkin least-squares problem $\Psi(X)K \approx \Phi(Y)$ between two spatial observable dictionaries.  The authors argue that when this finite section is rectangular the natural objects are singular modes rather than eigenmodes. Our construction is orthogonal to that program in two respects. First, PRONE is a discrete-time approximation of the Koopman operator itself, assembled from snapshot pairs $Y = F(X)$; WEDMD is a continuous-time approximation of the generator, and the weak form we employ is a temporal one, built on using compactly supported test functions $\psi_k(t)$ and 
integration by parts, whose purpose is precisely to avoid the
derivative estimation that a generator formulation would otherwise require. 

There is thus no dictionary substitution that reduces PRONE to WEDMD.  The dictionary $\Phi$ in PRONE consists of
functions of the pushed-forward state, whereas our $\psi_k$ are
functions of time. Second, the rectangularity of our pencil
$\lambda G\mathbf{w} = -D\mathbf{w}$ arises in the test-function
direction of a generalized eigenvalue problem, not from an asymmetric map between spatial observable spaces.  The eigenvalues of the pencil are the generator eigenvalues we seek, so eigen-analysis, rather than the singular-mode analysis appropriate to a rectangular spatial section, remains the correct tool. Finally, the noise robustness that
motivates WEDMD arises from temporal averaging over process and measurement noise, examined in Section~\ref{sec:linear}; PRONE does not treat stochastic data, and while WDMD also filters measurement noise, it does so for the linear state operator rather than the nonlinear generator setting considered here.

 With regard to establishing the range of applicability and accuracy of the WEDMD method, as a first check, we prove that our method produces exact results on linear systems.  Further, we examine the impact of noise by way of applying WEDMD on time series generated by linear stochastic ordinary differential equations (SDEs).  As observed, the weak formulation provides a natural averaging that readily mitigates the influence of noise.  

We then consider two nonlinear systems, the damped Duffing equation and the Van der Pol oscillator, and look at the performance of WEDMD in each system.  As demonstrated, our method does an excellent job of generating enough good approximations to Koopman modes so that WEDMD can even be used for forecasting purposes in noisy data.  Noise in this work consists of two classes.  The first and easier to address is additive Gaussian noise.  The second is Brownian motion, which as expected, does show the limits of WEDMD since SDEs can generate significant and nonlinear perturbations of the underlying dynamics.  For the Duffing equation, the attraction to a fixed point makes the impact of noise in the SDE case less significant, and good forecasts are still available.  On the other hand, the global nature of the limit cycle of the Van der Pol oscillator makes accurate prediction substantially more challenging and beyond the scope of the present work.  Addressing this shortcoming and developing further advances and analyses of the WEDMD method is the subject of future work.  

The structure of this paper is as follows.  In Section 2, we present the core of the WEDMD method.  This section also includes the analysis of linear problems.  In Section 3, we apply the WEDMD method to the two nonlinear systems described above.  A Discussion and Conclusion section then summarizes the results of the paper.  
\section{Weak-form Extended Dynamic Mode Decomposition}

Suppose we have some unknown dynamical system say
\[
\dot{{\bf y}} = {\bf f}({\bf y}), ~ {\bf y}(0)={\bf x}\in \mathbb{R}^{N_{s}},
\]
with the corresponding flow map ${\bf y}(t) = \varphi_{t}({\bf x})$.  In appropriately defined spaces \cite{LasotaMackey1994ChaosFractalsNoise}, the Koopman operator $\mathcal{K}^{t}$ where 
\[
\mathcal{K}^{t}g({\bf x}) = g\left(\varphi_{t}({\bf x}) \right),
\]
has the corresponding semigroup generator $\mathcal{L}$ with 
\[
\mathcal{L}g = {\bf f}({\bf x}) \cdot \nabla_{{\bf x}} g.
\]
The natural connection then between $\mathcal{K}^{t}$ and $\mathcal{L}$ comes by way of the eigenvectors and eigenvalues of $\mathcal{L}$ since 
\[
\mathcal{L}\phi = \lambda \phi
\]
immediately implies that 
\[
\mathcal{K}^{t}\phi({\bf x}) = e^{t\lambda}\phi({\bf x}).
\]
So, by finding the Koopman eigenfunctions via the generator, we can find the affiliated time-evolving Koopman operator.  

Of course, in the modern Koopman operator related literature we suppose that we have data ${\bf y}(t_{j})$ but not the corresponding vector field ${\bf f}({\bf x})$.  Following \cite{KLUS2020132416, Kaiser_2021}, one can try to still work with $\mathcal{L}$ by estimating $\dot{\bf y}_{j}$ and then solving some version of the regression problem
\[
\lambda \varphi({\bf y}_{j}) = \dot{{\bf y}}_{j} \cdot \nabla_{{\bf x}} \varphi({\bf y}_{j})
\]
by assuming that some dictionary of functions, say $\left\{{\bf f}_{l}^{(d)}({\bf x})\right\}_{l=1}^{N_{d}}$ well approximates the Koopman modes $\varphi({\bf x})$.  


Following the WSINDy ethos though, we want to avoid differentiation altogether by way of moving to weak forms of the problem.  Thus, we introduce the compactly supported test functions $\psi_{k}(t)$, $k=0,\cdots,K$, and integrate the eigenvalue problem so that 
\[
\lambda \int_{0}^{t_{f}} \varphi({\bf y}(t))\psi_{k}(t)dt = \int_{0}^{t_{f}}\psi_{k}(t)\left(\dot{{\bf y}} \cdot \nabla_{{\bf x}} \varphi({\bf y})\right)dt.
\]
Using
\begin{align*}
\int_{0}^{t_{f}}\psi_{k}(t)\left(\dot{{\bf y}} \cdot \nabla_{{\bf x}} \varphi({\bf y})\right)dt = & \int_{0}^{t_{f}}\psi_{k}(t)\frac{d}{dt}\varphi({\bf y}(t))dt\\
= & - \int_{0}^{t_{f}} \varphi({\bf y}(t)) \frac{d\psi_{k}}{dt} dt,
\end{align*}
if we then complement this with a dictionary assumption for $\varphi({\bf x})$ of the form 
\[
\varphi({\bf x}) \approx \sum_{l=1}^{N_{d}}w_{l}f^{(d)}_{l}({\bf x}) = \Theta({\bf x}){\bf w},
\]
we then get the system of equations
\[
\lambda \sum_{l=1}^{N_{d}}w_{l} \int_{0}^{t_{f}} f^{(d)}_{l}({\bf y}(t))\psi_{k}(t)dt = - \sum_{l=1}^{N_{d}}w_{l}\int_{0}^{t_{f}} f^{(d)}_{l}({\bf y}(t)) \frac{d\psi_{k}}{dt} dt.
\]

Specifying that we have $K+1$ compactly supported test functions $\psi_{k}(t)$ with supports contained within $(0,t_{f})$, we then get the $(K+1)\times N_{d}$ generalized eigenvalue problem 
\[
\lambda {\bf G} {\bf w} = - {\bf D} {\bf w}.
\]
where
\[
{\bf G} = \int_{0}^{t_{f}} \boldsymbol{\Psi}(t)\boldsymbol{\Theta}({\bf y}(t))dt \approx \delta t\sum_{j=1}^{N_{T}-1}\boldsymbol{\Psi}_{j}\boldsymbol{\Theta}({\bf y}_{j})
\]
and
\[
{\bf D} = \int_{0}^{t_{f}} \dot{\boldsymbol{\Psi}}(t)\boldsymbol{\Theta}({\bf y}(t))dt \approx \delta t\sum_{j=1}^{N_{T}-1}\dot{\boldsymbol{\Psi}}_{j}\boldsymbol{\Theta}({\bf y}_{j})
\]
where 
\[
\boldsymbol{\Psi}(t) = \left(\psi_{0}(t)\cdots \psi_{K}(t) \right)^{T}.
\]
Further, if we have discrete data $\left\{{\bf y}_{j}\right\}_{j=0}^{N_{T}}$ sampled uniformly with step size $\delta t$, we can control the support size of each function $\psi_{k}(t)$ so that discretization via the Trapezoid Rule keeps quadrature errors to nominal levels.  

With regard to practical implementations, since we must ultimately discretize the integrals, we define 
\[
{\bf D} = \boldsymbol{\dot{\Psi}}_{d}\boldsymbol{\Theta}_{d}, ~ {\bf G} = \boldsymbol{\Psi}_{d}\boldsymbol{\Theta}_{d}
\]
where the $d$ subscript denotes that quadrature via the Trapezoid Rule has been applied so that given data set $\left\{{\bf y}_{j}\right\}_{j=0}^{N_{T}}$ we have 
\[
\boldsymbol{\Psi}_{d} = \begin{pmatrix}
\psi_{1}(t_{0}) & \cdots & \psi_{1}(t_{f}) \\ 
\vdots & \ddots & \vdots \\ 
\psi_{K}(t_{0}) & \cdots & \psi_{K}(t_{f}) 
\end{pmatrix}, ~ \boldsymbol{\Theta}_{d} = \begin{pmatrix}
f_{1}^{(d)}({\bf y}_{0}) & \cdots & f_{N_{d}}^{(d)}({\bf y
}_{0}) \\ 
\vdots & \ddots & \vdots \\ 
f_{1}^{(d)}({\bf y}_{N_{T}}) & \cdots & f_{N_{d}}^{(d)}({\bf y}_{N_{T}}) 
\end{pmatrix}. 
\]
Forming and then solving this discretized generalized eigenvalue problem forms the core of what we call the WEDMD method.

\subsection{Test-Function Selection and Approximation Theory}
To understand the advantage of this approach, assuming that we have equispaced data $\left\{{\bf y}_{j}\right\}_{j=0}^{N_{T}}$ sampled at the rate $\delta t$, we choose $\psi_{k}(t)$ such that 
\[
\psi_{k}(t)=\tilde{\psi}\left(\frac{2(t-c_{k})}{r}\right), 
\]
where following prior WSINDy approaches, standard bump functions of the form 
\[
\tilde{\psi}(t) = \left\{ \begin{array}{rl} (1-t^{2})^{p} &  t\in[-1,1] \\ 0 & t\notin [-1,1]\end{array}\right.
\]
with $p>1$ are used.  The overlap and centers are determined relative to a choice of radius $r$, overlap parameter $0<s<1$, and polynomial power $p$ with 
\[
\psi_{k}(t_{\ast}) = \psi_{k+1}(t_{\ast}) = s, 
\]
so that $t_{\ast} = (c_{k}+c_{k+1})/2$ and 
\[
c_{k+1} - c_{k} = \frac{r}{2}\sqrt{1-s^{1/p}}.
\]
We naturally require that
\[
[0,t_{f}] \subset \cup_{k=0}^{K}[c_{k}-r/2,c_{k}+r/2], ~ c_{0}=r/2, ~ c_{K} = t_{f}-r/2,
\]
so that 
\[
K = \left\lceil \frac{2(t_{f}-r)}{r\sqrt{1-s^{1/p}}} \right\rceil.
\]
Moreover, we observe that for fixed $r$ that as $s\rightarrow 1^{-}$ then $K\rightarrow \infty$.  Due to our use of numerical quadrature, we finally require that each interval $[c_{k}-r/2,c_{k}+r/2]$ have at least one discrete time $t_{j}$ in the interior.  

Clearly however, our set of test functions cannot exhaust $L_{2}([0,t_{f}])$ since $K$ is finite.  To that end, if we define the subspace
\[
T(s) = \left\{\sum_{k=0}^{K}\left<f,\psi_{k}\right> \psi_{k}, ~ f\in L_{2}([0,t_{f}])\right\}
\]
then naturally we have 
\[
L_{2}([0,t_{f}]) = T(s) \oplus T(s)^{\perp}.
\]
Defining the affiliated operator ${\bf F}_{s}:L_{2}([0,t_{f}])\rightarrow T(s)$ such that
\[
{\bf F}_{s}f = \sum_{k=0}^{K}\left<f, \psi_{k}\right>\psi_{k},
\]
then 
\[
\left<{\bf F}_{s}(f),f\right> = \sum_{k=0}^{K}\left|\left<f, \psi_{k}\right>\right|^{2}
\]
and thus we find that ${\bf F}_{s}$ is a bounded, positive semi-definite operator.  To find a bound on ${\bf F}_{s}$, for $f \in L_{2}([0, t_{f}])$, if we look at the quantity 
\[
\left<{\bf F}_{s}f,f\right> = \sum_{k=0}^{K}\left|\left<f, \psi_{k}(t)\right>\right|^{2} 
\]
by Minkowski's inequality, we have 
\[
\left(\sum_{k=0}^{K}\left|\left<f, \psi_{k}(t)\right>\right|^{2}\right)^{1/2} \leq \int_{0}^{t_{f}} dt|f(t)|\left(\sum_{k=0}^{K}\psi^{2}_{k}(t)\right)^{1/2},
\]
so that using Cauchy-Schwarz we find 
\[
\sum_{k=0}^{K}\left|\left<f,\psi_{k}(t)\right>\right|^{2} \leq \gnorm{f}_{1}^{2}\leq C_{K}\gnorm{f}_{2}^{2},
\]
where 
\begin{align*}
C_{K} = & \sum_{k=0}^{K}\int_{0}^{t_{f}}\psi_{k}^{2}(t)dt, \\
= & \frac{r}{2}(K+1)\int_{-1}^{1}(1-t^{2})^{2p}dt.  
\end{align*}
Therefore we have 
\[
\sup_{\gnorm{f}_{2}=1}\left<{\bf F}_{s}f,f\right> \leq C_{K}.
\]
We can of course make this bound uniform in $K$ by using $\psi_{k}(t)/\sqrt{K+1}$ as test functions instead.  

As to what we cannot represent in $T(s)^{\perp}$, from the positive semi-definiteness, we see that 
\[
\text{ker}({\bf F}_{s}) = \left\{g\in L_{2}([0,t_{f}]), \left<g,\psi_{k}\right>=0, ~ k=0,\cdots,K\right\}.
\]
Clearly $\text{ker}({\bf F}_{s})\subset T(s)^{\perp}$, and for $g\in T(s)^{\perp}$ we have     
\[
\left<{\bf F}_{s}(g),g\right> = 0,
\]
so $g\in \text{ker}({\bf F}_{s})$.  Thus, for a bounded self-adjoint operator, we have the typical orthogonal decomposition
\[
L_{2}([0,t_{f}]) = \text{rng}({\bf F}_{s}) \oplus \text{ker}({
\bf F
}_{s}).
\]
Further, considering the $s\rightarrow 1^{-}$ limit, for $g(t)\in \text{ker}({\bf F}_{s})$, from 
\[
\left<g,\psi_{k}\right> = 0, ~ k=0,\cdots, K,
\]
using Parseval's equality, we then find 
\[
\left<g,\psi_{k}\right> = \frac{r}{2\sqrt{t_{f}}}\sum_{m-\infty}^{\infty}\hat{g}_{m}C_{m} e^{2\pi i m c_{k}/t_{f}},  
\]
where 
\[
C_{m} = \int_{-1}^{1}\cos\left(\frac{\pi m r \tau}{t_{f}}\right)(1-\tau^{2})^{p}d\tau.
\]
For $p \geq 1$, by Lebesgue's Lemma, $\hat{g}_{m}C_{m}$ is an absolutely convergent series if $\hat{g}_{m}\in l_{2}$.  Therefore $\bar{\psi}(t)$ given by
\[
\bar{\psi}(t) = \frac{r}{2\sqrt{t_{f}}}\sum_{m-\infty}^{\infty}\hat{g}_{m}C_{m} e^{2\pi i m t/t_{f}}
\]
is a continuous (increasingly smooth for larger $p$) $t_{f}$-periodic function.  As $s\rightarrow 1^{-}$ or $K\rightarrow \infty$, then $\bar{\psi}(t)$ would have to have an infinite number of zeros on $[0,t_{f}]$, which is only possible if $\hat{g}_{m}C_{m}=0$ for all $m\in \mathbb{Z}$.  

Thus, as $s\rightarrow 1^{-}$, the question of the size and content of $\text{ker}({\bf F}_{s})$ hinges on when $C_{m}=0$. Writing a proof at this point is challenging, but we can explore the issue through a few examples.  If $p=1$, we find that $C_{m}=0$ only for those $m\neq 0$ such that  
\[
\tan\left(\frac{\pi m r}{t_{f}}\right) = \frac{\pi m r}{t_{f}}.
\]
Thus, aside from non-generic choices of $r$ and $t_{f}$, this should not be satisfied for any $m\neq 0$.  Looked at from another vantage point, if say $r/t_{f}=1/M$, $M\in \mathbb{N}$, $M>1$, then $C_{m}\neq 0$ for $|m| < M/2$.  Thus, for small supports relative to the entire interval, $\text{ker}({\bf F}_{s})$ would only consist of relatively high frequency signals.    

Ultimately, these issues matter with regard to the properties of ${\bf G}$.  Treating the continuous time case, we observe that if there were a vector ${\bf w}$ such that ${\bf G}{\bf w}={\bf 0}$ then
\[
{\bf G}{\bf w} = \begin{pmatrix}\left<\psi_{0},\boldsymbol{\Theta}({\bf y}){\bf w}\right> \\ \vdots \\ \left<\psi_{K},\boldsymbol{\Theta}({\bf y}){\bf w}\right>\end{pmatrix} = {\bf 0}.
\]
Following the arguments above, it should be possible to choose $r$, $s$, and $p$ to ensure ${\bf G}$ has full rank.  This allows us to pass to the more straightforward eigenvalue problem in the form
\[
{\bf G}^{-P}{\bf D}{\bf w} = -\lambda {\bf w},
\]
where ${\bf G}^{-P}$ is the Moore-Penrose pseudoinverse of ${\bf D}$.  Finally, our approximate Koopman eigenfunctions are given by $\varphi_{l}({\bf x}) = \boldsymbol{\Theta}({\bf x}){\bf w}_{l}$.

\subsection{Residual Error Analysis}

Of course, there is the question of the initial assumption that we can write Koopman eigenfunctions in terms of a finite number of dictionary functions.  Using our approximation to the Koopman eigenfunctions, our closure approximation can be written as 
\[
\mathcal{K}^{t}\boldsymbol{\Theta}({\bf x}){\bf W} = \boldsymbol{\Theta}({\bf x}){\bf W}e^{t\boldsymbol{\Lambda}} + {\bf R}(t,{\bf x}).
\]
Differentiating and integrating against test function $\psi_{k}(t)$ and then performing integration-by parts as above, we get 
\begin{multline*}
-\int_{0}^{t_{f}}\dot{\psi}_{k}(t)\boldsymbol{\Theta}({\bf y}(t)){\bf W} dt = \int_{0}^{t_{f}}\psi_{k}(t)\boldsymbol{\Theta}({\bf y}(t)){\bf W}\boldsymbol{\Lambda}  dt \\- \int_{0}^{t_{f}}{\bf R}(t,{\bf x})\left(\boldsymbol{\Lambda}\psi_{k}(t) + \dot{\psi}_{k}(t)\right)dt
\end{multline*}

Our method then determines ${\bf W}$ and the corresponding eigenvalues $\boldsymbol{\Lambda}$ so that 
\begin{align*}
\int_{0}^{t_{f}}{\bf R}(t,{\bf x})\left(\boldsymbol{\Lambda}\psi_{k}(t) + \dot{\psi}_{k}(t)\right)dt = 0, ~ k=0, \cdots, K.  
\end{align*}
Per our choice of $\psi_{k}(t)$, this can be rewritten as 
\[
\int_{-1}^{1}d\tau ~{\bf R}\left(c_{k} + \frac{r}{2}\tau, {\bf x}\right)\left(\boldsymbol{\Lambda} \tilde{\psi} - \dot{\tilde{\psi}}\right) = 0, ~ k=0, \cdots, K.
\]
So as in \cite{MessengerBortz2021WSINDy}, we see that our method is a Galerkin approximation that chooses ${\bf W}$ and $\boldsymbol{\Lambda}$ so as to make the error orthogonal to the linear combination of test functions $\lambda \psi_{k} + \dot{\psi}_{k}$.  

Further, we find that as $r\rightarrow 0^{+}$, by dominated convergence we should have  
\begin{align*}
\lim_{r\rightarrow 0^{+}}\int_{-1}^{1}d\tau ~{\bf R}\left(c_{k} + \frac{r}{2}\tau,{\bf x}\right)\left(\boldsymbol{\Lambda} \tilde{\psi} - \dot{\tilde{\psi}}\right) = & 0\\
{\bf R}\left(t,{
\bf x
} \right)\boldsymbol{\Lambda}\int_{0}^{1}\tilde{\psi}(\tau)d\tau = & 0.  
\end{align*}
Thus, by making the test-function supports arbitrarily small while letting the centers $c_{k}$ pass to a continuum limit, for fixed initial condition ${\bf x}$, we can essentially eliminate the residual in our closure approximation.  While an infinite resolution data result, we see the advantage of a weak formulation insofar as it directly lends itself towards promoting convergence for a particular initial condition ${\bf x}$.  On the other hand, the degree to which this result generalizes for other values of ${\bf x}$ is still limited by the degree to which the observables form a closed or nearly closed subspace of the Koopman operator.   




\subsection{Analysis of Linear Systems}\label{sec:linear}
That the underlying eigenvalue problem of Koopman analysis leads to an approximate eigenvalue problem is not surprising, but it does introduce new issues not necessarily encountered in the original WSINDy method.  It is therefore instructive to look at a linear system of the form 
\[
\dot{{\bf y}} = {\bf A}{\bf y}, ~ {\bf y}(0) = {\bf x}.  
\]
We have ${\bf y}(t)  = e^{{\bf A}t}{\bf x}$, and if we only include the canonical observation functions, i.e. $f_{l}({\bf x})=x_{l}$, then our method generates the system
\[
{\bf G} = \int_{0}^{t_{f}}\boldsymbol{\Psi}(t) \left(e^{{\bf At}}{\bf x}\right)^{T}dt, ~ {\bf D} = \int_{0}^{t_{f}}\dot{\boldsymbol{\Psi}}(t) \left(e^{{\bf At}}{\bf x}\right)^{T}dt
\]

If we suppose that ${\bf A} = {\bf V}\boldsymbol{\Lambda}{\bf V}^{-1}$, then our pencil problem
\[
{\bf G}{\bf W}\boldsymbol{\Lambda} = - {\bf D}{\bf W}
\]
can be rewritten as 
\[
\tilde{{\bf G}}\tilde{{\bf W}}\boldsymbol{\Lambda} = - \tilde{{\bf D}}\tilde{{\bf W}},
\]
with 
\[
\tilde{{\bf G}} = \int_{0}^{t_{f}}\boldsymbol{\Psi}(t)\tilde{{\bf x}}^{T}e^{\boldsymbol{\Lambda}t}dt, ~ \tilde{{\bf D}} = \int_{0}^{t_{f}}\dot{\boldsymbol{\Psi}}(t)\tilde{{\bf x}}^{T}e^{\boldsymbol{\Lambda}t}dt
\]
with $\tilde{{\bf x}} = {\bf V}^{-1}{\bf x}$. In turn then, we find that  
\begin{align*}
\tilde{{\bf G}}_{kn} = & \tilde{x}_{n}\int_{0}^{t_{f}}\psi_{k}(t)e^{\lambda_{n}t}dt,\\
= & \frac{r\tilde{x}_{n}}{2}e^{\lambda_{n}c_{k}}\int_{-1}^{1}\tilde{\psi}(t)e^{r\lambda_{n}t/2}dt,
\end{align*}
and
\begin{align*}
\tilde{{\bf D}}_{kn} = & \tilde{x}_{j}\int_{0}^{t_{f}}\dot{\psi}_{k}(t)e^{\lambda_{n}t}dt,\\
= &  \tilde{x}_{n} e^{\lambda_{j}c_{k}}\int_{-1}^{1}\dot{\tilde{\psi}}(t)e^{r\lambda_{n}t/2}dt.
\end{align*}
We find immediately then that
\[
\tilde{{\bf D}}_{kn} = -\lambda_{n}\tilde{{\bf G}}_{kn}.  
\]
Over a linear system then, the only source of error in our method would be either from quadrature error or noise in the data.  All other things being equal though, we should get the exact Koopman operator $e^{{\bf A}t}$.    

More interesting then is how the weak approach addresses the presence of noise.  If we suppose that we sample data from the linear stochastic equation
\[
d{\bf y} = {\bf A}{\bf y}dt + {\bf C}d{\bf W}_{t}, ~ {\bf y}(0)={\bf x}
\]
where ${\bf C}$ is a real-valued full rank square matrix, ${\bf W}_{t}$ is standard Brownian motion, and we treat ${\bf x}$ as deterministic.  Then using standard arguments from the Ito calculus, we have that 
\[
{\bf G} = \int_{0}^{t_{f}}\boldsymbol{\Psi}(t)\left(e^{{\bf A}t}{\bf x}\right)^{T}dt + \int_{0}^{t_{f}}\boldsymbol{\Psi}(t)\left(\int_{0}^{t}e^{{\bf A}(t-s)}{\bf C}d{\bf W}_{s}\right)^{T}dt, 
\]
and thus 
\[
\text{Var}({\bf G}_{kn}) = r^{2}\int_{0}^{1}\int_{0}^{1}\tilde{\psi}(t_{1})\tilde{\psi}(t_{2})\mathbb{E}\left[\tilde{{\bf y}}^{T}(c_{k}+rt_{1}/2)\tilde{{\bf y}}(c_{n}+rt_{2}/2) \right]dt_{1}dt_{2},
\]
where $\tilde{{\bf y}}(t) = {\bf y}(t) - e^{{\bf A}t}{\bf x}$ and
\[
\mathbb{E}\left[\tilde{{\bf y}}^{T}(t_{1})\tilde{{\bf y}}(t_{2}) \right] = \int_{0}^{\text{min}(t_{1},t_{2})}\text{tr}\left(e^{{\bf A}(t_{1} - s)}{\bf C}{\bf C}^{T} e^{{\bf A}^{T}(t_{2} - s)}\right)ds.
\]
Likewise, we can show that 
\[
\text{Var}({\bf D}_{kn}) = \int_{-1}^{1}\int_{-1}^{1}\dot{\tilde{\psi}}(t_{1})\dot{\tilde{\psi}}(t_{2})\mathbb{E}\left[\tilde{{\bf y}}^{T}(c_{k}+rt_{1}/2)\tilde{{\bf y}}(c_{n}+rt_{2}/2) \right]dt_{1}dt_{2}.
\]
Therefore, for the linear observables considered here, integrating
against test functions applies a de-facto low-pass filter over the
noise: the variances of ${\bf G}$ and ${\bf D}$ depend only on filtered
correlations of the data, in contrast to direct estimation of
$\dot{{\bf y}}$, which unless very carefully done tends to amplify
variance even for linear problems.

We emphasize that this filtering interpretation is specific to the
linear-observable setting. The identity used to pass from
$\psi_{k}\,\dot{{\bf y}}\cdot\nabla_{{\bf x}}\varphi$ to a boundary term
is the ordinary chain rule; for an It\^o process it must be replaced by
It\^o's lemma, which contributes an additional drift
$\tfrac{1}{2}\nabla_{{\bf x}}\cdot\left({\bf C}{\bf C}^{T}\nabla_{{\bf x}}
\varphi\right)$. For the canonical observables $f_{l}({\bf x}) = x_{l}$
this term vanishes identically, since $\nabla_{{\bf x}}^{2}\varphi = 0$,
which is precisely why process noise enters above only through the
zero-mean fluctuation $\tilde{{\bf y}}$ that the averaging suppresses.
For a nonlinear dictionary this is no longer the case: the It\^o
correction is a deterministic shift of the operator being estimated,
not a zero-mean perturbation, and it is not removed by integrating
against test functions. In that regime the weak form consistently
estimates the generator of the \emph{stochastic} dynamics --- the
Kolmogorov backward generator --- rather than that of the underlying
deterministic flow. We return to this distinction when interpreting the
stochastic Van der Pol results in Section~\ref{sec:vanderpol}.

\subsection{Determining Test Function Support}

Per our convergence analysis from above, and directly in line with the existing WSINDy literature, the choice of test function support, again denoted by $r$, is a critical hyperparameter that has an outsize influence on the overall performance of our proposed approach.  Beyond its role in improving the convergence of the Galerkin approximation, the support width, as explored in \cite{MessengerBortz2021_WSINDyPDE}, serves as a means for managing noise through the implicit use of the test functions as low-pass filters.  

Moreover, there is the question as to how well the approximated Koopman eigenfunctions and values function in a time-stepping capacity.  We can track this error via the function $\mathcal{E}_{c}(r)$ defined as
\[
\mathcal{E}^{2}_{c}(r) = \text{max}_{1\leq l\leq N_{d}}\frac{\sum_{j=1}^{N_{T}}\left|\boldsymbol{\Theta}({\bf y}_{j}){\bf w}_{l}(r) - \boldsymbol{\Theta}({\bf x}){\bf w}_{l}(r)e^{j\delta t \lambda_{l}(r)}\right|^{2}}{\sum_{j=1}^{N_{T}}\left|\boldsymbol{\Theta}({\bf y}_{j}){\bf w}_{l}(r)\right|^{2}}.
\]
$\mathcal{E}_{c}(r)$ measures the relative accuracy to which a given eigenvector and eigenvalue act as Koopman eigenfunctions and values over the length the time series.  We then determine $r$ by minimizing $\mathcal{E}_{c}(r)$.  To limit our search, we first find an estimate for $r$ using the changepoint detection approach outlined in \cite{MessengerBortz2021WSINDy}.  

\subsection{Test Function Selection and Complete Algorithm}

We now discuss the choice of test function.  Assuming that $r$ and $s$ have been chosen, following \cite{MessengerBortz2021_WSINDyPDE}, we can add further specificity to the method by introducing the parameter $0<\tau_{l}\ll 1$ which is the value of the test function one step back from the edge of its support.  Defining $m_s = \lfloor r/(2\delta t)\rfloor$, we choose $p$ so that 
\[
\left(1-\left(1 - \frac{1}{m_{s}}\right)^{2}\right)^{p} = \tau_{l}.  
\]

\begin{algorithm}[H]
\caption{WEDMD Algorithm}
\label{wedmd}
\noindent INPUT: $\delta t$, $\left\{{\bf y}_{j}\right\}_{j=0}^{N_{T}}$, $s\in(0,1)$, $0<\tau_{l} \ll 1$, $\left\{f_{l}^{(d)}({\bf  x})\right\}_{l=1}^{N_{D}}$
\begin{algorithmic}[1]
	\Procedure{Generate ${\bf G}$, ${\bf D}$}{}	
        \State INITIALIZE: Find $r$ and $p$ via the Changepoint Detection 
        Method \cite{MessengerBortz2021_WSINDyPDE}. 
        \State COMPUTE: $r \leftarrow \text{arg min}_{\tilde{r}} \mathcal{E}_{c}(\tilde{r})$.
        \State COMPUTE: $p \leftarrow \log(\tau_{l})/\log((2m_{s}-1)/m_{s}^{2})$, ~ $m_{s}=\lfloor r/2\delta t\rfloor$.
        \State BUILD: $\boldsymbol{\Psi}_{d}$, $\boldsymbol{\dot{\Psi}}_{d}$, $\boldsymbol{\Theta}_{d}$.
	\State RETURN: ${\bf G}$, ${\bf D}$
    \EndProcedure
\end{algorithmic}
\end{algorithm}

\section{Examples}

Throughout the examples, the underlying non-noisy data is generated via fourth order Runge-Kutta.  For the SDE data, that is generated via a second order stochastic Runge-Kutta scheme.  Throughout, we fix the observable dictionary to be products of monomials with power up to $N_{ord}$.  In the two dimensional examples that follow then, $N_{d} = (N_{ord}+1)^{2}$ and 
\[
f_{l}^{(d)}({\bf x}) = x_{1}^{l_{1}}x_{2}^{l_{2}}, ~ l = (N_{ord}+1)l_{1} + l_{2}, ~ 0\leq l_{j} \leq N_{ord}.  
\]
Results using Legendre and Hermite polynomials were also generated, but no significant difference in performance was observed.  
\subsection{The Damped Duffing Equation}

\begin{figure}
\centering
\includegraphics[width=.75\textwidth]
{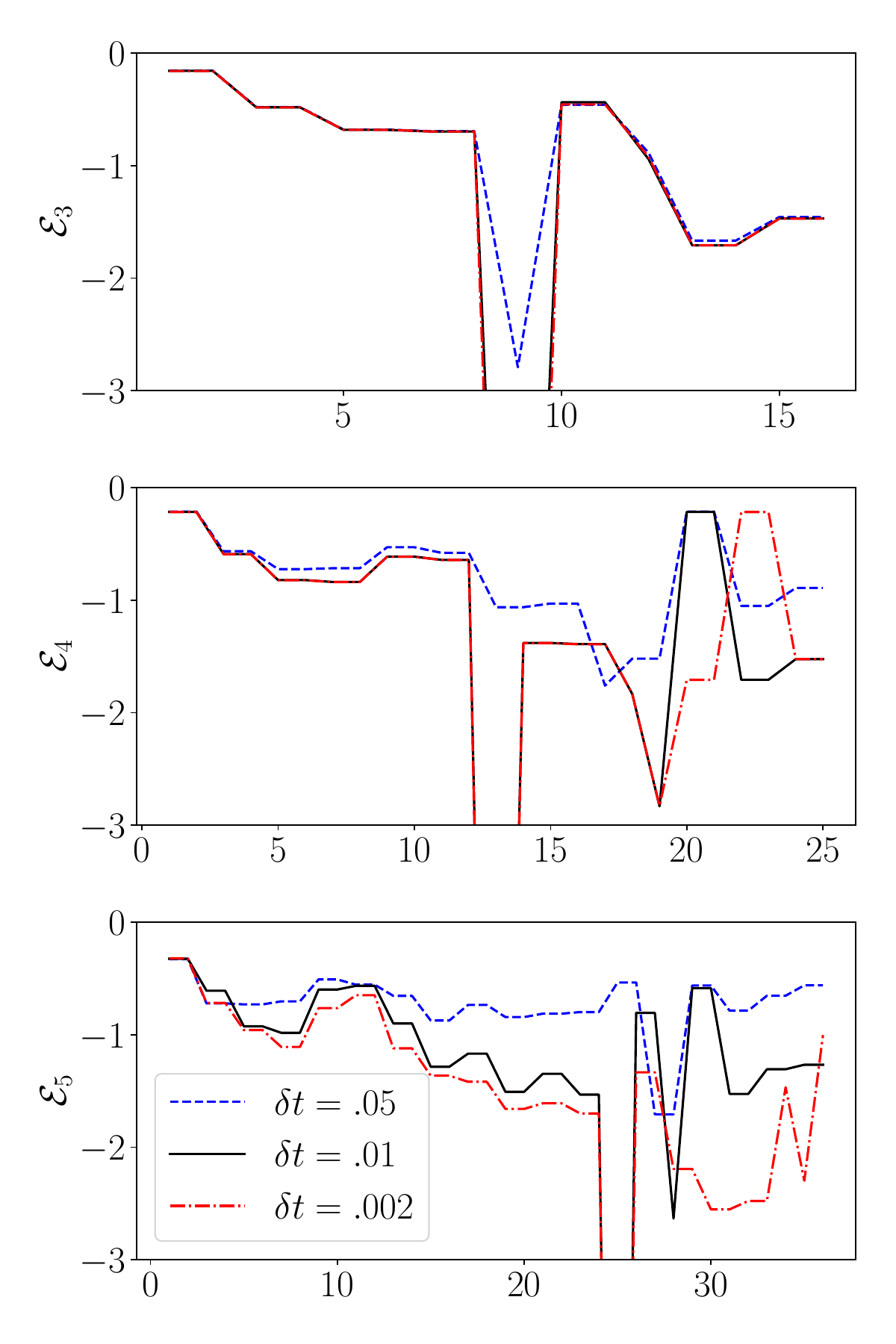}
\caption{Mode error for $N_{ord}=3 ~(\text{top}), 4~(\text{middle}), 5~(\text{bottom})$ and $\delta t = .05, .01, .002$.}
\label{fig:duffing_prob_mode_error}
\end{figure}

We begin our analysis on the damped Duffing equation given by
\[
\dot{y}_{1} = y_{2}, ~ \dot{y}_{2} = -y_{1} + \frac{y_{1}^{3}}{6} - \mu y_{2}.
\]
Choosing weak damping with $\mu=.1$, we look at data for $0\leq t \leq 20$, which for the relatively weak damping, is long enough to see sufficiently interesting dynamics.  In Figure \ref{fig:duffing_prob_mode_error}, we look at varying $N_{ord}$ from $N_{ord}=3$ to $N_{ord}=5$ and also how varying $\delta t$ influences mode error measured by $\mathcal{E}_{N_{ord}}$, which is defined the same way as $\mathcal{E}_{c}(r)$ but now with the emphasis on varying $N_{ord}$.  On clean data, the method shows tendencies towards convergence, with larger proportions of high accuracy modes appearing as we increase $N_{ord}$.  Nevertheless, we also see that the method struggles to generate uniform levels of performance, reflecting the limitations of picking arbitrary spaces of observables not otherwise adapted to the particular data set.    

\begin{figure}
\centering
\includegraphics[width=.75\textwidth]
{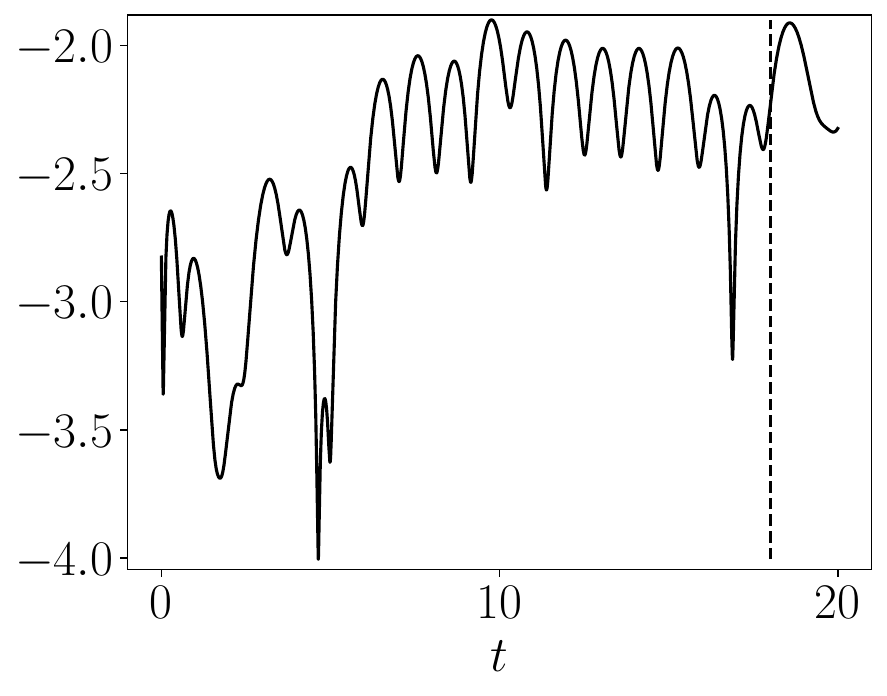}
\caption{Relative error using a reduced mode model for ${\bf y}(t)$ in the Duffing equation with $N_{ord}=5$, $\delta t = .002$, and the filtering over those modes with $\mathcal{E}_{5}<.1$.  The vertical bar indicates the beginning of model prediction.}
\label{fig:duffing_prob_prediction}
\end{figure}

Moreover, while the method does not compute uniformly accurate Koopman modes, by simply filtering over sufficiently accurate modes, we can get excellent models of the dynamics.  Specifically, choosing $N_{ord}=5$ and $\delta t = .002$, we select those modes $({\bf w}_{l},\lambda_{l})$ such that $\mathcal{E}_{5,l} < .1$, 
so that we only select modes with at worst 10 \% error.  From this reduced set of modes, we can then generate a time-evolving model of the dynamics via the approximation
\[
{\bf y}^{T}(t) \approx \boldsymbol{\Theta}({\bf x}){\bf W}_{r}e^{\boldsymbol{\Lambda}_{r}t}{\bf K},
\]
where ${\bf K}$ is found through standard regression via the Frobenius norm.  To test this, we generate modes using data for $0\leq t\leq 18$ and then use our reduced mode model to predict dynamics from $18\leq t \leq 20$.  We see the results of this in Figure \ref{fig:duffing_prob_prediction}.  As seen, the relative error from our model relative to the numerically generated solution to the Duffing equation is at worst 1 \%, even in the prediction regime.  So, even though we cannot compute uniformly accurate modes, simply filtering over reasonably good ones produces outstanding results for a nonlinear problem.  

To examine the impact of noise, fixing $N_{ord}=5$ and $\delta t = .01$, we look at two noise models.  The first is additive noise, in which the data ${\bf y}_{j}$ is given by 
\[
{\bf y}_{j} = {\bf y}^{(tr)}_{j} + \sigma \boldsymbol{\epsilon}_{j}, ~ \boldsymbol{\epsilon}_{jn} \sim \mathcal{N}(0,1).
\]
The second model adds Brownian motion on to the damped Duffing equation so that we work with data generated by the stochastic differential equation (SDE)
\[
dy_{1} = y_{2}dt + \sigma dW^{(1)}_{t}, ~ dy_{2} = \left(-y_{1} + \frac{y_{1}^{3}}{6} - \frac{y_{2}}{10}\right)dt + \sigma dW^{(2)}_{t}.
\]
In our case, we suppose the noise is uncorrelated across dimensions with the variance again given by the parameter $\sigma$.  We use the second-order stochastic Runge-Kutta method to simulate paths perturbed by Brownian motion. 

\begin{figure}
\begin{tabular}{cc}
\includegraphics[width=.5\textwidth]{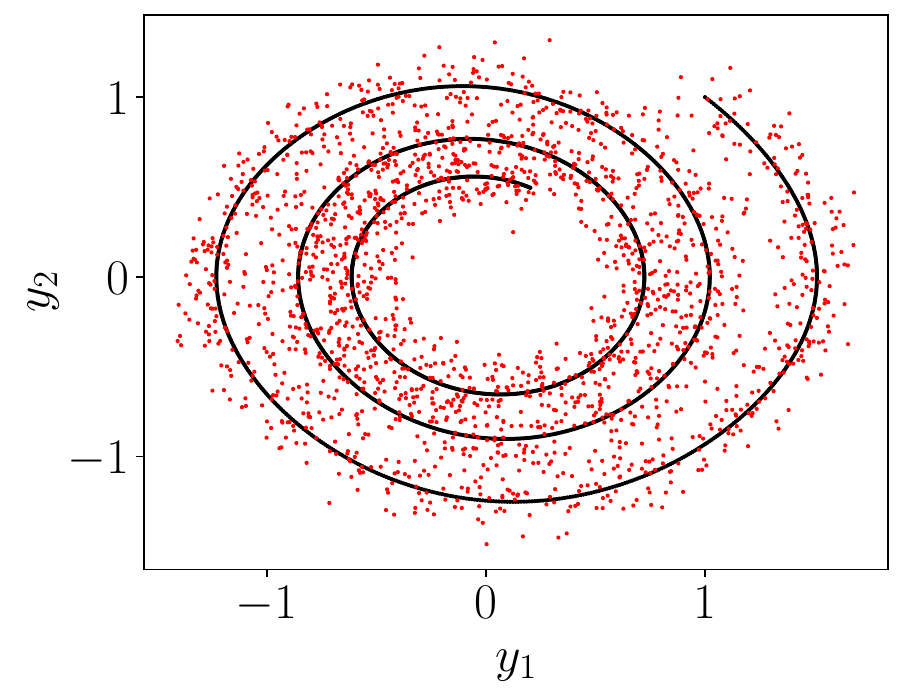} & \includegraphics[width=.5\textwidth]{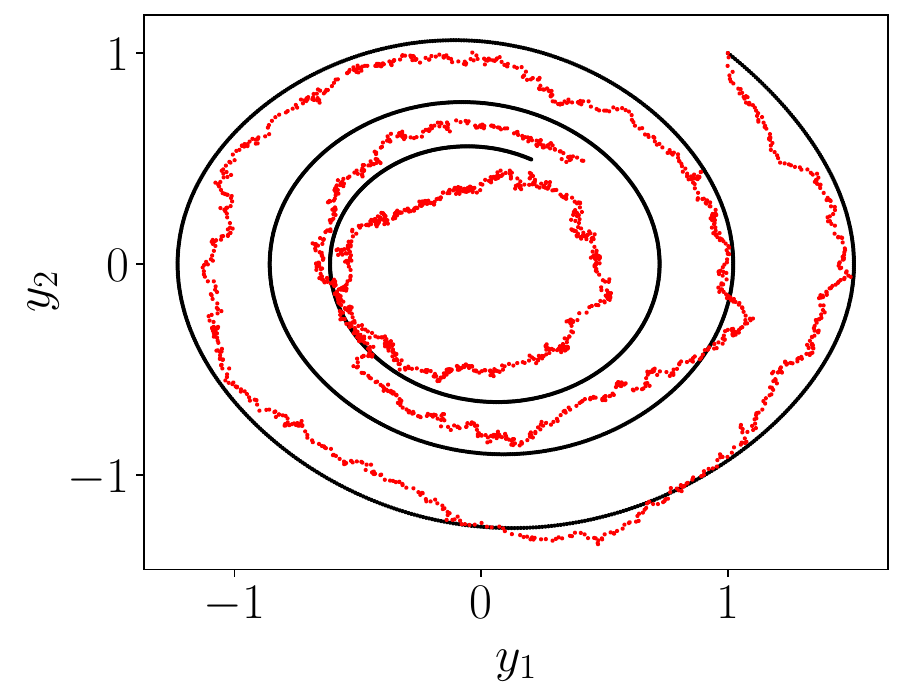}\\
(a) & (b) \\ 
\includegraphics[width=.5\textwidth]{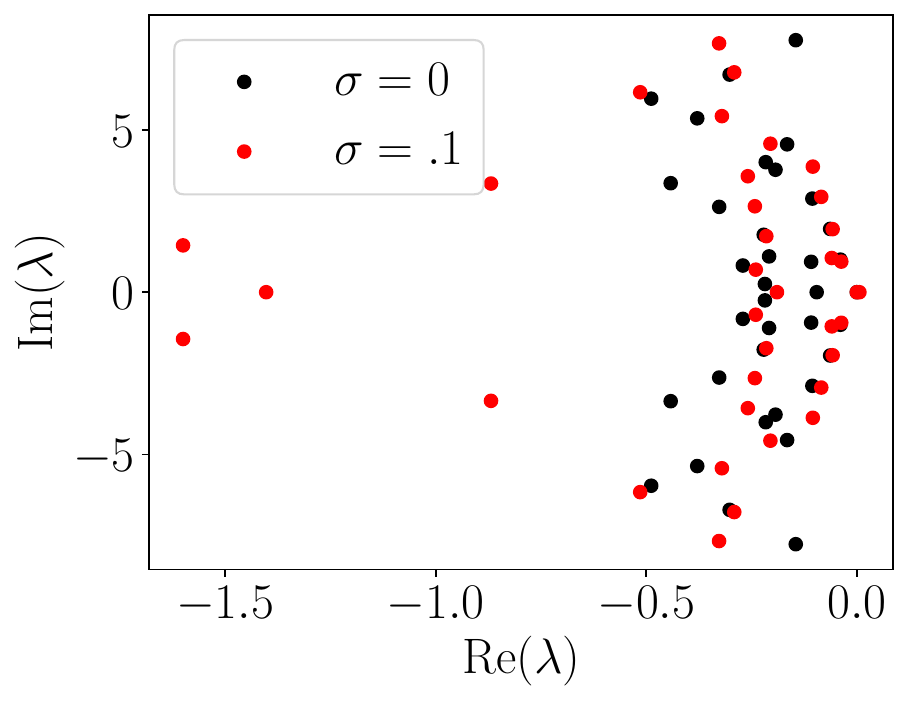} & \includegraphics[width=.5\textwidth]{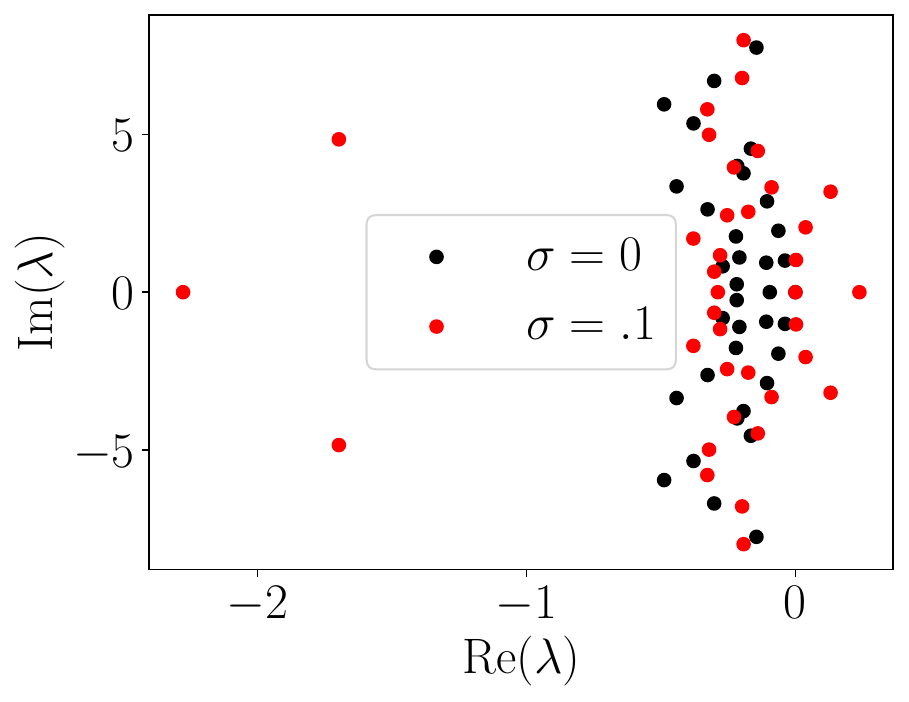}\\
(c) & (d)\\
\includegraphics[width=.5\textwidth]{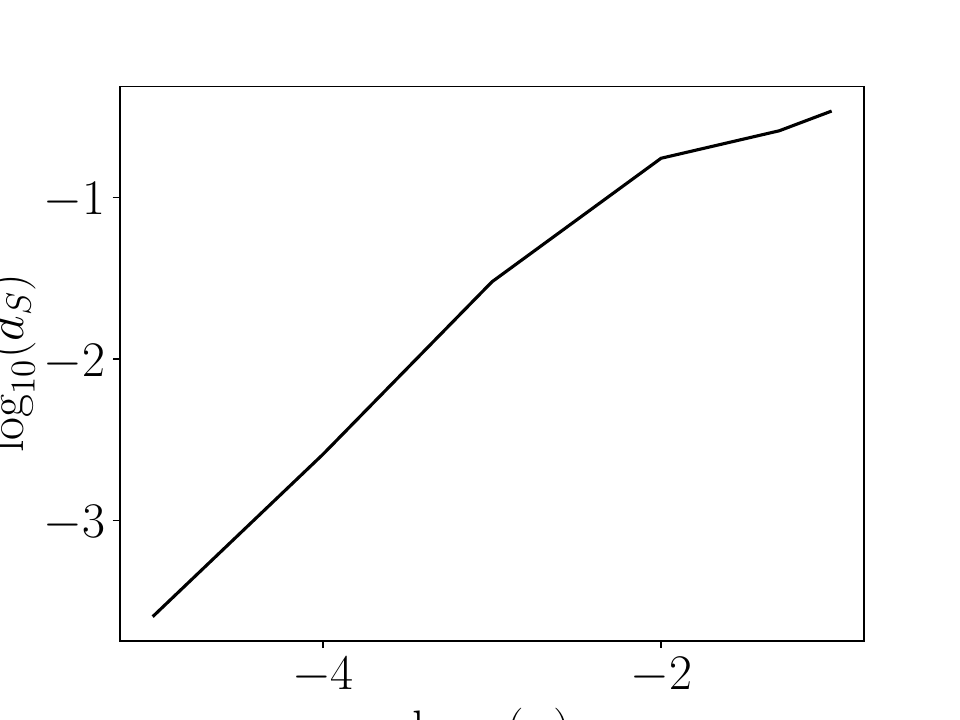} & \includegraphics[width=.5\textwidth]{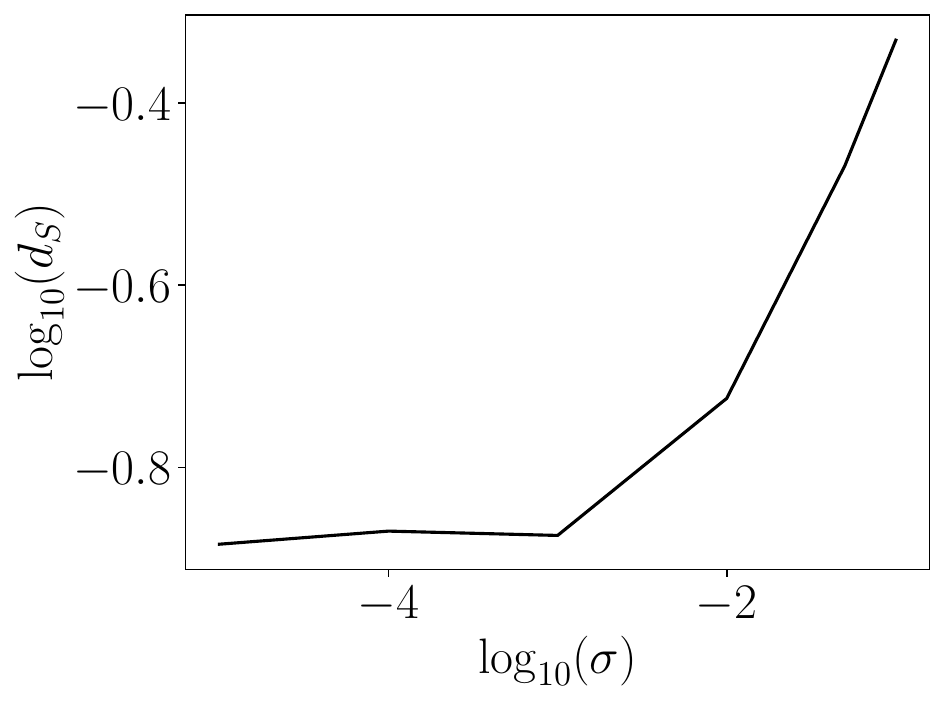}\\
(e) & (f)
\end{tabular}
\caption{For the Duffing equation, noisy data generated by additive noise (left panels) and an SDE (right panels).  Typical realizations with $\sigma=.1$ are shown in (a) and (b) with corresponding computed spectra in (c) and (d).  Average distances between clean and noisy spectra are shown in (e) and (f).}
\label{fig:duffing_spectral_distance}
\end{figure}

We can then look at how varying $\sigma$ causes deviation in the computed spectra.  To compute this, we define the distance metric $d_{S}$ where 
\[
d_{S} = d\left(\sigma_{ns},\sigma_{cl}\right) = \max_{1\leq l \leq N_{d}}\min_{1\leq n \leq N_{d}}\left|\lambda_{l,ns} - \lambda_{n,cl}\right| 
\]
where $\sigma_{ns}$ is the noisy spectrum and $\sigma_{cl}$ the clean spectrum.  Figure \ref{fig:duffing_spectral_distance} shows typical dynamics and the distance between the average separation in the spectra from the no-noise case, i.e. $\sigma=0$, and the case with nontrivial noise.  Average distances are computed across ten realizations.  As seen, the method is able to handle noise smoothly in both cases, though the SDE generates stronger perturbations of the spectra and thus overall larger distances $d_{S}$.  As can be seen comparing spectra in Figure \ref{fig:duffing_spectral_distance} (c) and (d), much of the separation comes from clear outliers so that questionable spectra could be identified using basic k-means or other group classification schemes.  

\begin{figure}[!h]
\centering
\begin{tabular}{cc}
\includegraphics[width=.5\textwidth]{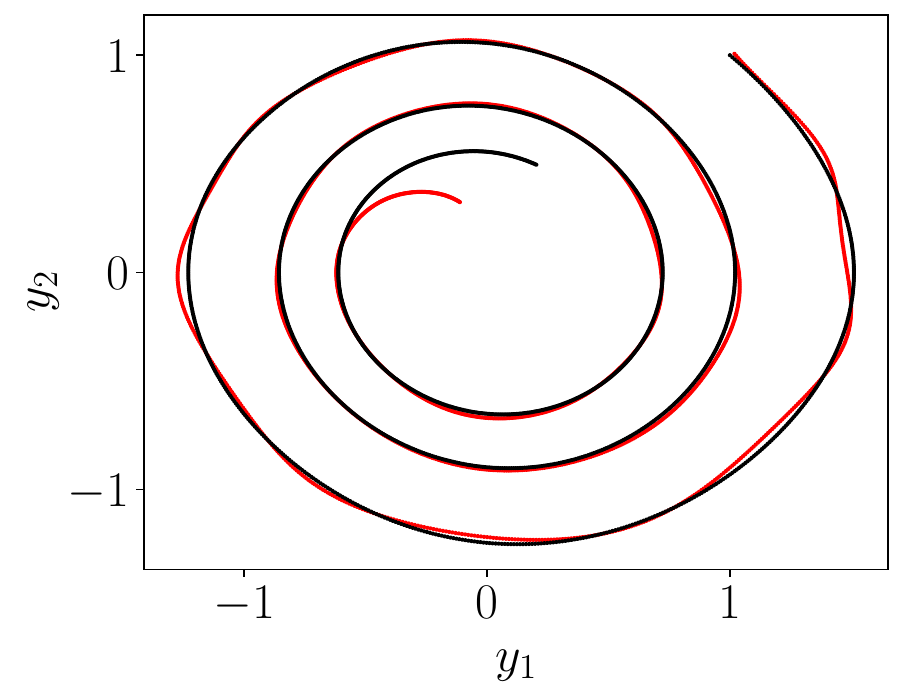} & \includegraphics[width=.5\textwidth]{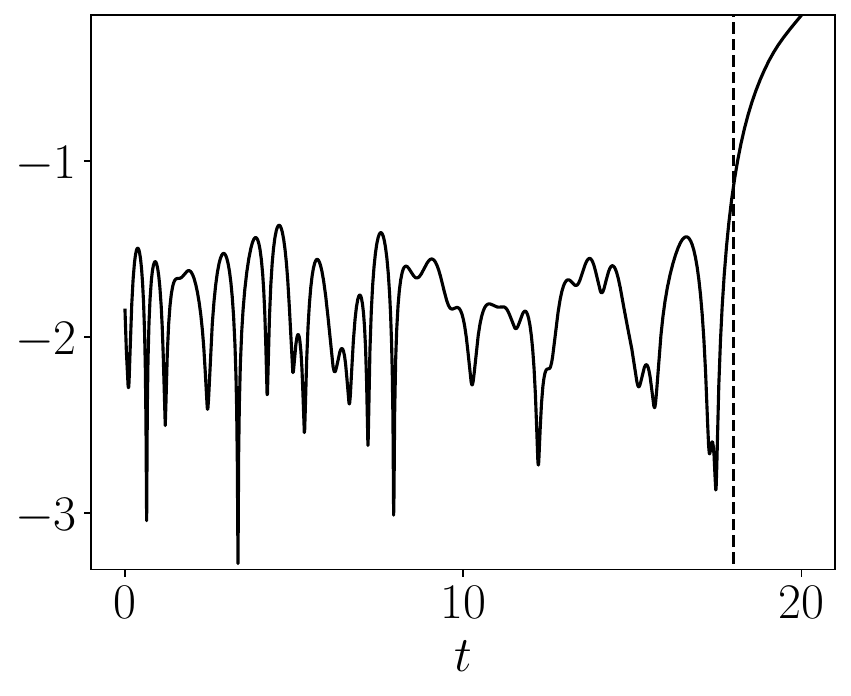}
\end{tabular}
\caption{Comparison of Duffing dynamics generated by reduced model using additive noise perturbed data (red) to noise-free data (black) (left) and error in reduced-order model (right).  Noise level is $\sigma=.1$.  Vertical line in right figure indicates beginning of model prediction.}
\label{fig:duffing_add_noise_prediction}
\end{figure}

Looking at the prediction problem for the additive noise case with $\sigma=.1$ and $N_{ord}=5$, if we just filter over modes such that the error $\mathcal{E}_{l}$ satisfies $\mathcal{E}_{l}<1.$, then we get the results in Figure \ref{fig:duffing_add_noise_prediction}.  As we see, the averaging properties of the Galerkin method allow for surprising levels of accuracy, making our method at the least an excellent additive noise filtration scheme.  

\begin{figure}[!h]
\centering
\begin{tabular}{cc}
\includegraphics[width=.5\textwidth]{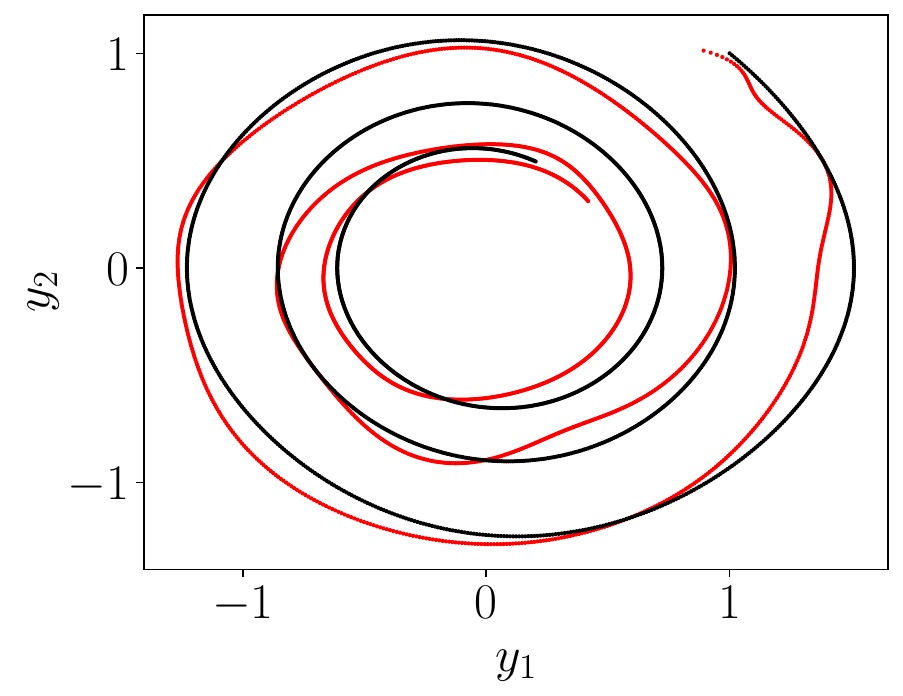} & \includegraphics[width=.5\textwidth]{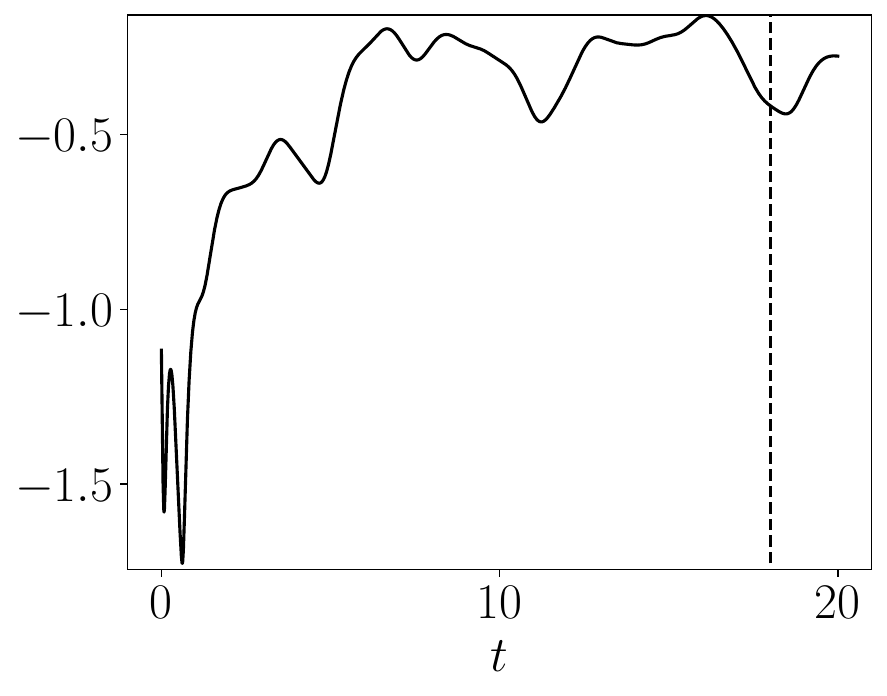}
\end{tabular}
\caption{Comparison of Duffing dynamics generated by reduced model using SDE noise perturbed data (red) to noise-free data (black) (left) and error in reduced-order model (right).  Noise level is $\sigma=.1$.  Vertical line in right figure indicates beginning of model prediction.}
\label{fig:duffing_sde_noise_prediction}
\end{figure}

For the SDE case, the results are not as good, though given the degree of perturbation, the results are still satisfactory.  Again letting $\sigma=.1$, we see in Figure \ref{fig:duffing_sde_noise_prediction} that if we increase $N_{ord}=6$ and filter modes so that $\mathcal{E}_{l}<2$, then we get some reasonable approximation of the true dynamics from the SDE data. 

\subsection{Van der Pol Oscillator}\label{sec:vanderpol}

\begin{figure}
\centering
\includegraphics[width=.75\textwidth]
{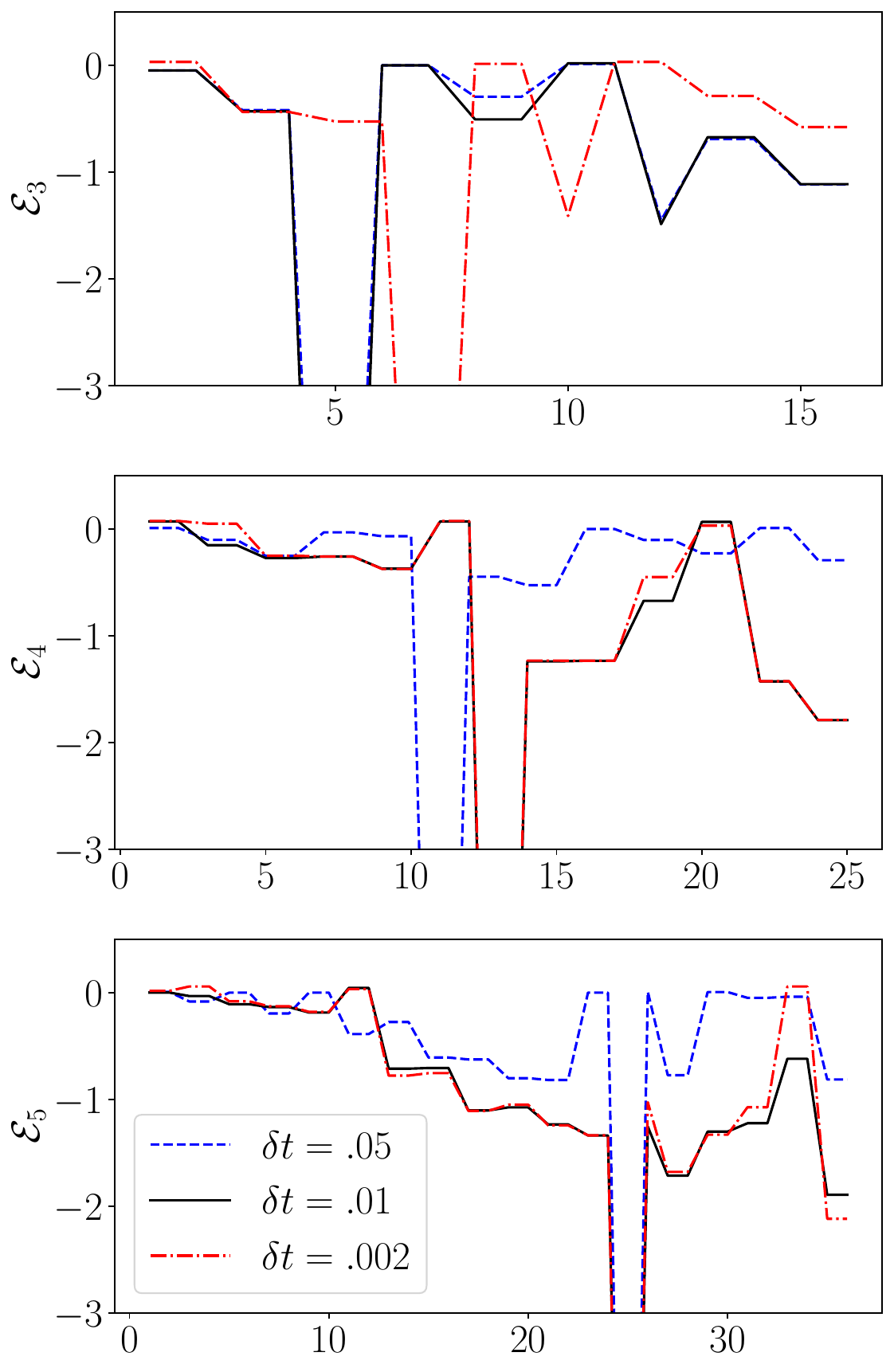}
\caption{Mode error for the Van der Pol oscillator for $N_{ord}=3 ~(\text{top}), 4~(\text{middle}), 5~(\text{bottom})$ and $\delta t = .05, .01, .002$.}
\label{fig:van_der_pol_prob_mode_error}
\end{figure}

\begin{figure}
\centering
\includegraphics[width=.75\textwidth]
{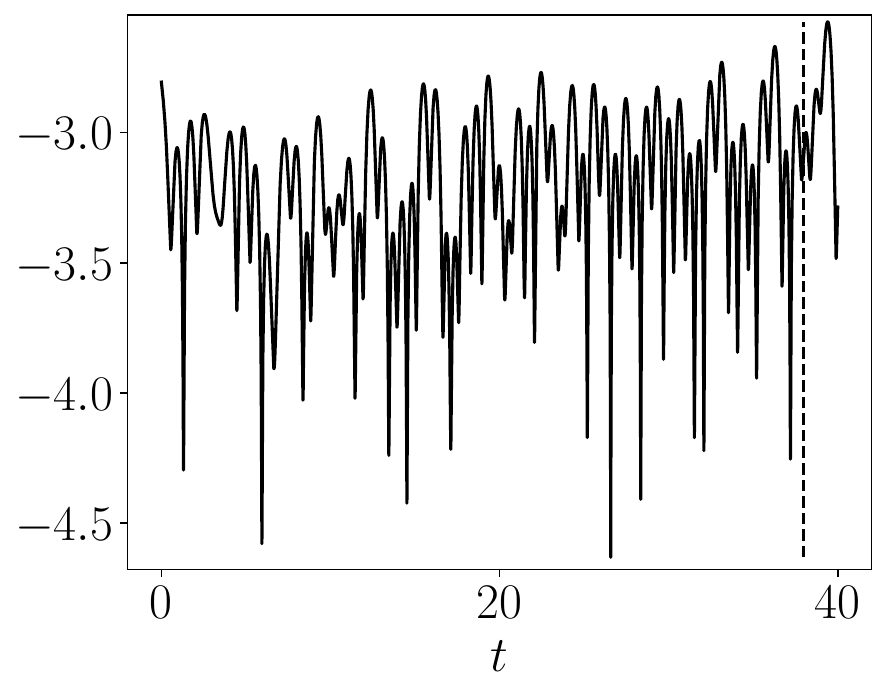}
\caption{Relative error using a reduced mode model for ${\bf y}(t)$ in the Van der Pol equation with $N_{ord}=5$, $\delta t = .002$, and the filtering over those modes with $\mathcal{E}_{5}<.1$.  The vertical bar indicates the beginning of model prediction.}
\label{fig:van_der_pol_prob_prediction}
\end{figure}

We now look at the Van der Pol oscillator given by the two-dimensional system
\[
\dot{y}_{1} = y_{2}, ~ \dot{y}_{2} = \mu(1-y_{1}^{2})y_{2}-y_{1}. \]
Avoiding the more singular case, we let $\mu=.1$ with data being generated from $0\leq t \leq 40$ so that enough time passes for the limit cycle to attract the dynamics.  The distinction here with the prior example is that the long time limit is a periodic orbit instead of just a single attracting point.  Referring to Figure \ref{fig:van_der_pol_prob_mode_error}, this added complexity shows up by way of higher mode errors than seen for the Duffing equation example.  Nevertheless, the WEDMD approach is able to compute a large number of relatively excellent modes, and again, this improves as $N_{ord}$ is increased and $\delta t$ decreases, indicating some degree of convergence.  
Repeating the mode reduction experiments from above, we see in Figure \ref{fig:van_der_pol_prob_prediction}, that we get even better performance.  Our method then is able to find excellent approximations of Koopman modes so that with only a minimal selection criteria, we can generate good forecasting models.    

\begin{figure}
\begin{tabular}{cc}
\includegraphics[width=.5\textwidth]{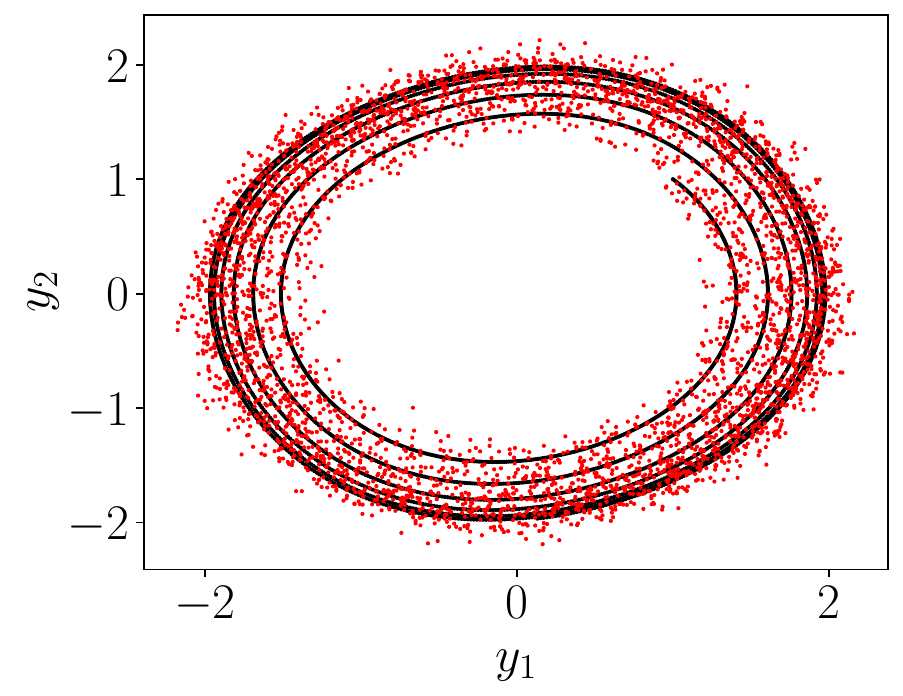} & \includegraphics[width=.5\textwidth]{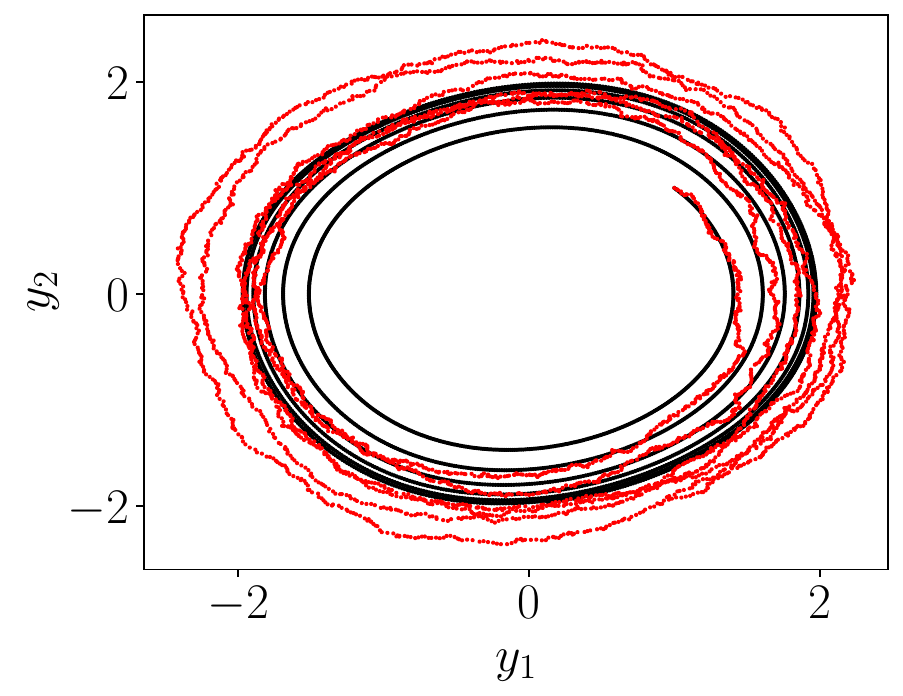}\\
(a) & (b) \\ 
\includegraphics[width=.5\textwidth]{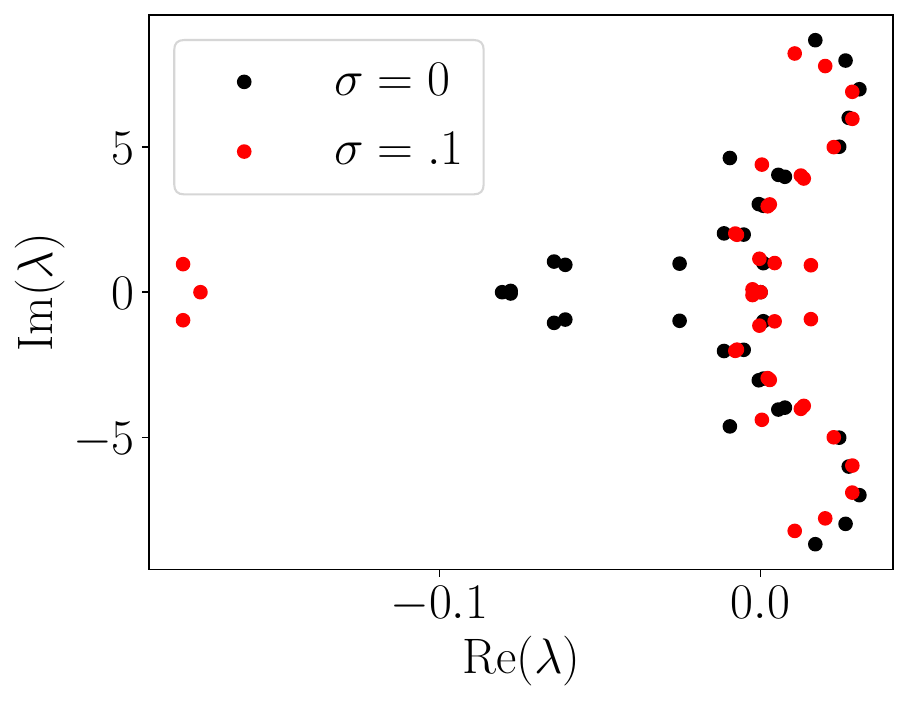} & \includegraphics[width=.5\textwidth]{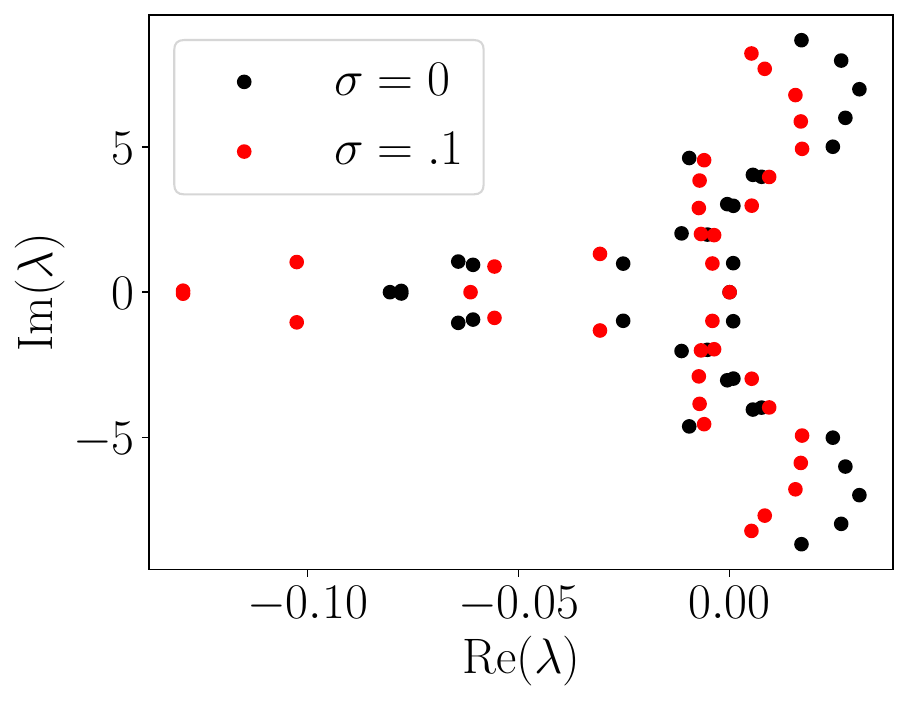}\\
(c) & (d)\\
\includegraphics[width=.5\textwidth]{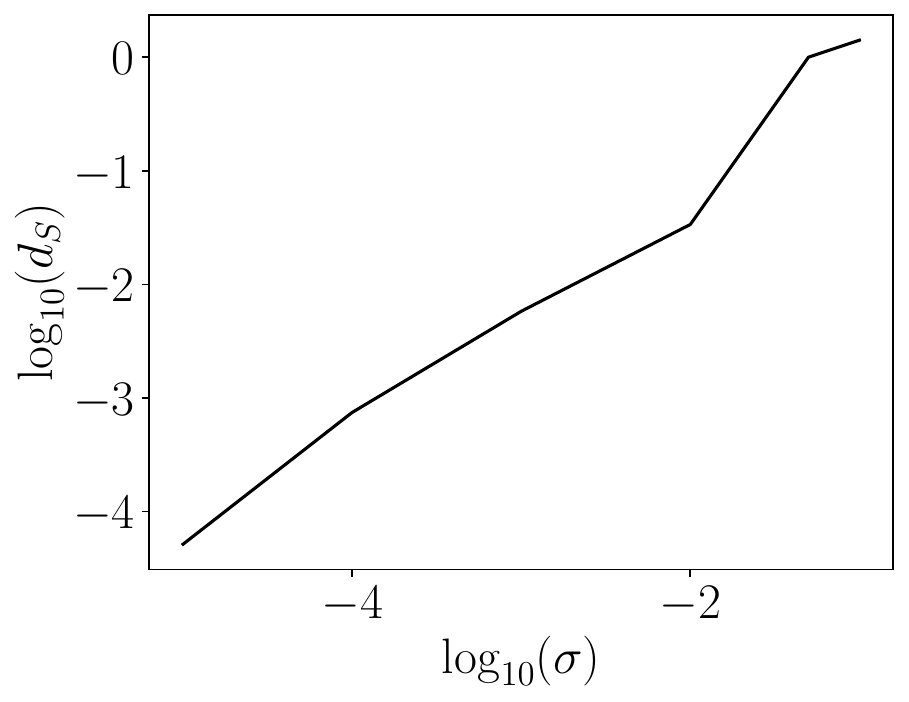} & \includegraphics[width=.5\textwidth]{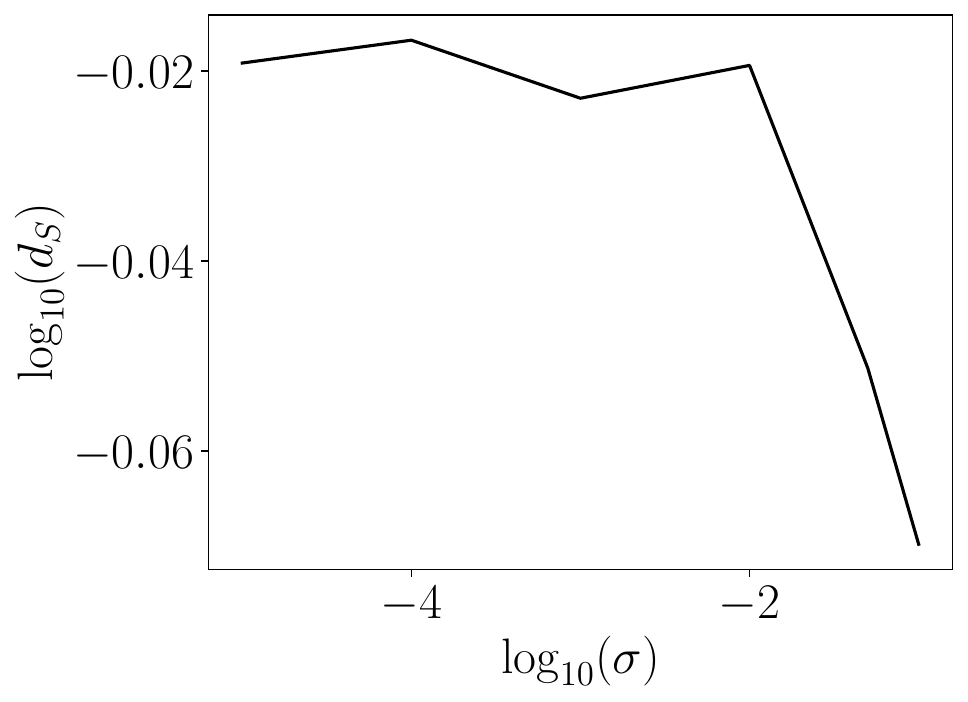}\\
(e) & (f)
\end{tabular}
\caption{For the Van der Pol equation, noisy data generated by additive noise (left panels) and an SDE (right panels).  Typical realizations with $\sigma=.1$ are shown in (a) and (b) with corresponding computed spectra in (c) and (d).  Average distances between clean and noisy spectra are shown in (e) and (f).}
\label{fig:van_der_pol_spectral_distance}
\end{figure}

Adding noise to the problem shows that the global limit cycle creates sharp distinctions in how the method performs on different types of noise.  Referring to Figure \ref{fig:van_der_pol_spectral_distance}, in the additive case, the averaging property of the weak approach makes the computation of spectra and its overall response to noise virtually identical in performance to the Duffing equation.  Errors are somewhat more pronounced and larger, however, k-means would again readily isolate outliers in the noisy spectra.  However, for the SDE case, the limit cycle makes computing spectra even for very small levels of noise very delicate with much more diffuse and difficult to isolate spreads in the computed spectra.  

\begin{figure}
\centering
\begin{tabular}{cc}
\includegraphics[width=.5\textwidth]{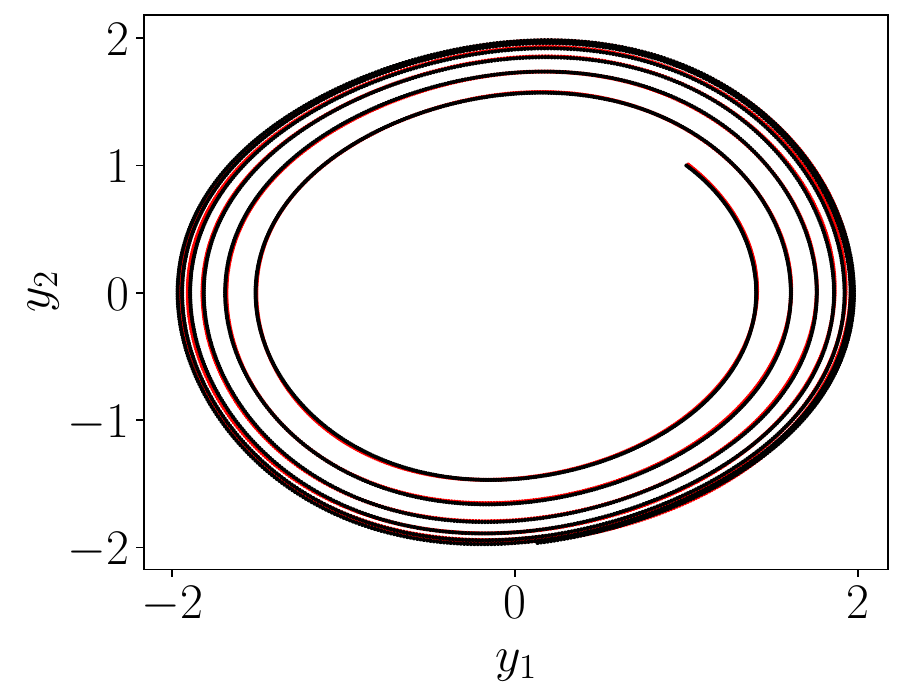} & \includegraphics[width=.5\textwidth]{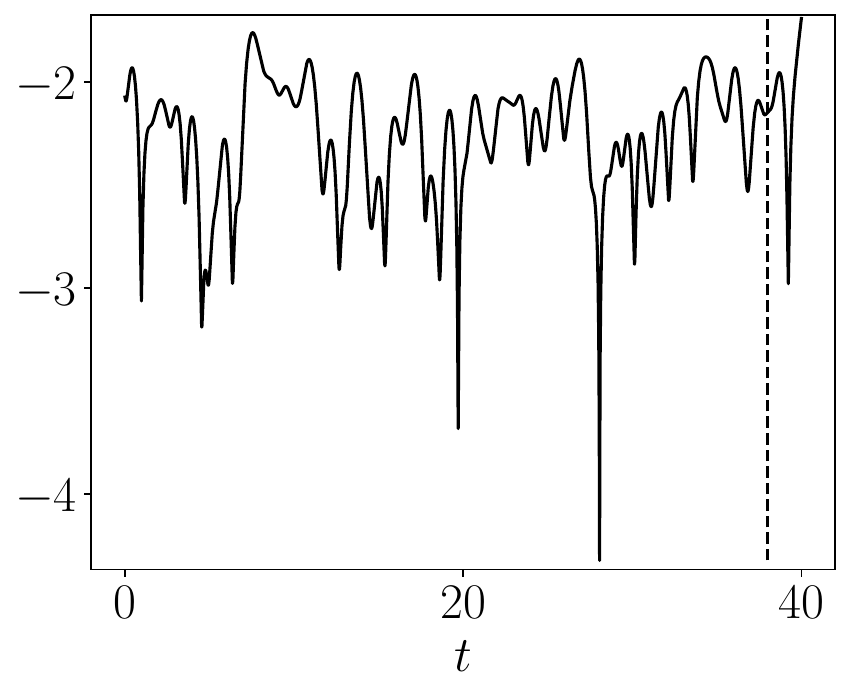}
\end{tabular}
\caption{Comparison of Van der Pol dynamics generated by reduced model using additive noise perturbed data (red) to noise-free data (black) (left) and error in reduced-order model (right).  Noise level is $\sigma=.1$.  Vertical line in right figure indicates beginning of model prediction.}
\label{fig:van_der_pol_add_noise_prediction}
\end{figure}

\begin{figure}
\centering
\includegraphics[width=.5\textwidth]{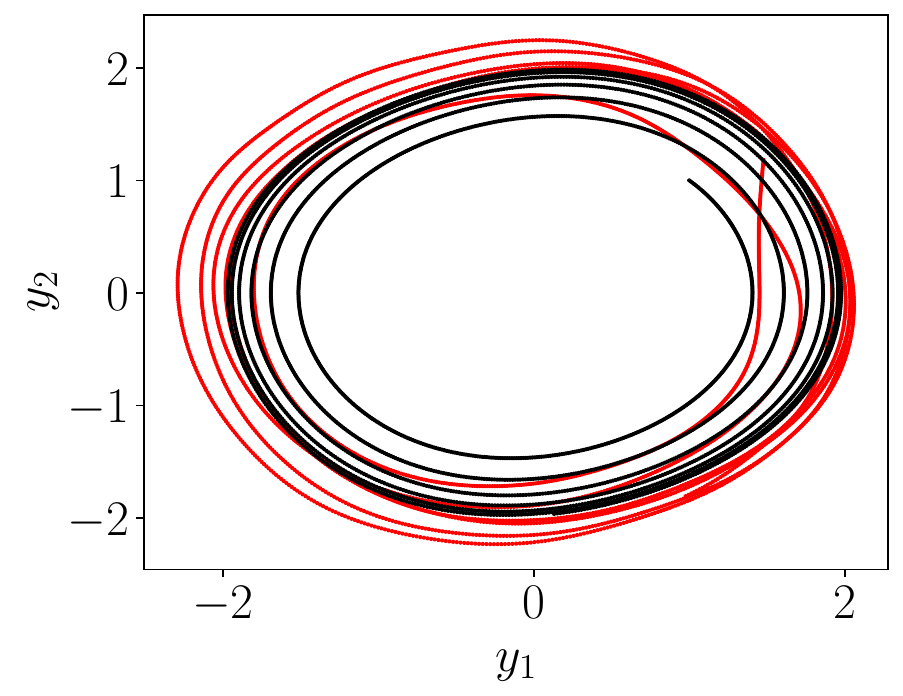} \caption{Comparison of Van der Pol dynamics generated by reduced model using SDE noise perturbed data (red) to noise-free data (black).  Noise level is $\sigma=.1$.}
\label{fig:van_der_pol_sde_noise_prediction}
\end{figure}

Looking at the prediction problem for the additive noise case with $\sigma=.1$ and $N_{ord}=5$, by filtering over those modes satisfying $\mathcal{E}_{l}<1.$, we get the results in Figure \ref{fig:van_der_pol_add_noise_prediction}.  As for the Duffing equation, the averaging properties of the Galerkin method produce exceptional levels of accuracy.  
However, in the SDE case, as seen in Figure \ref{fig:van_der_pol_sde_noise_prediction}, the method tracks the data not the underlying dynamics.  This reflects our poor approximation of the full generator of the stochastic dynamics as discussed at the end of Section 2.3.  

\section{Discussion and Conclusion}
In this work, we developed the Weak-form Extended Dynamic Mode Decomposition (WEDMD) method, a weak-form extension of EDMD inspired by recent advances in weak-form data-driven model discovery. By replacing pointwise evaluations with Galerkin projections onto compactly supported test functions, WEDMD avoids explicit differentiation and naturally incorporates averaging that mitigates the effects of noise. Our analysis showed that the method recovers the exact Koopman spectrum for linear deterministic systems up to quadrature error, while the weak formulation provides a systematic means of controlling approximation error through the choice of test-function support. Furthermore, the stochastic analysis demonstrated that the resulting operators depend on filtered correlations of the data, providing an inherent low-pass filtering effect absent from pointwise approaches.

The numerical examples confirmed these theoretical advantages. For both the damped Duffing equation and the Van der Pol oscillator, WEDMD generated increasingly accurate Koopman modes as the observable dictionary was enriched and the sampling interval was refined. Although not all computed modes exhibited uniformly low error, simple filtering based on reconstruction accuracy consistently identified subsets of high-quality modes that produced excellent reduced-order models with  correspondingly excellent forecasting performance. In particular, the method proved remarkably robust to additive Gaussian noise, maintaining accurate spectral estimates and predictive capability even when the data were substantially perturbed.

We also see limitations in the method, especially with regard to the impact of noise arising from SDEs. While WEDMD remained effective for stochastic Duffing dynamics, for the Van der Pol oscillator, Brownian motion produced significantly larger perturbations than additive measurement noise and reduced forecast accuracy. In this case, the global structure of the attracting limit cycle caused stochastic perturbations to generate diffuse Koopman spectra and forecasts that tracked the noisy trajectories rather than the underlying deterministic dynamics. Future work will focus on extending WEDMD to stochastic Koopman frameworks, where a more careful examination of the impacts of process noise and the concomitant measurement noise are addressed and better accounted for in our methodology.  A key future approach in this direction is developing adaptive observable selection strategies by way of autoencoders or other machine learning approaches. Nevertheless, the results presented here demonstrate that WEDMD provides a robust and flexible framework for Koopman spectral analysis and forecasting in the presence of noisy data.
\section*{Acknowledgments}
This work was supported in part by the NIGMS Division of Biophysics, Biomedical Technology and Computational Biosciences under grant R35GM149335 and in part by the DOE CHaRMNET Mathematical Multifaceted Integrated Capability Center (MMICC)   under grant DE-SC0023346.
\bibliography{wdmd_bibliography_siam}
\bibliographystyle{siam}

\end{document}